\documentclass[10pt,reqno]{article}

\newif\iffigures\figurestrue
\figuresfalse

\usepackage[colorlinks=false,pdfborderstyle={/W 1}]{hyperref}
\newif\ifhyper\IfFileExists{hyperref.sty}{\hypertrue}{\hyperfalse}

\ifhyper\usepackage{hyperref}
\def\hitem#1#2{\item[\hypertarget{#1}{#2}]\expandafter\gdef\csname LBL#1ITM\endcsname{#2}}
\def\iref#1{\hyperlink{#1}{\csname LBL#1ITM\endcsname}}
\else
\def\hitem#1#2{\item[{#2}]\expandafter\gdef\csname LBL#1ITM\endcsname{#2}}
\def\iref#1{{\csname LBL#1ITM\endcsname}}
\fi

\usepackage{amssymb,amsmath,amsthm,amsfonts,mathrsfs,ulem, cancel}

\usepackage{graphicx}
\usepackage[utf8]{inputenc}
\usepackage{enumitem}

\def\text#1{\hbox{#1}}

\def\N{{\mathbb N}}
\def\T{{\mathbb T}}

\def\eps{\epsilon}

\def\1{\mathbf{1}}

\newtheorem{theorem}{Theorem}
\newtheorem{lemma}[theorem]{Lemma}

\newtheorem{remark}{{\bf Remark}}
\newtheorem*{note}{\underline{Note}}

\long\def\note#1/{\ifdraft{\marginpar{{$\Longleftarrow$}} \bf [#1] }\fi}

\newcounter{my}
\stepcounter{my}

\def\1{1\!\! 1}

\def\T{T}

\theoremstyle{definition}
\numberwithin{equation}{section}
\numberwithin{figure}{section}

\title{A unimodular random graph with large upper growth and no growth}

\author{Péter Mester, Ádám Timár} 

\date{ \today}

\begin{document}

\maketitle

\begin{abstract}
We construct a unimodular random rooted graph with maximal degree $d\geq 3$ and upper growth rate $d-1$, which does not have a growth rate. Abért, Fraczyk and Hayes
showed that for a unimodular random tree, if the upper growth rate is at least $\sqrt{d-1}$, then the growth rate exists, and asked with some scepticism if this may hold for more general graphs. Our construction shows that the answer is negative. We also provide a non-hyperfinite example of a unimodular random graph with no growth rate. This may be of interest in light of a conjecture of Abért that unimodular Riemannian surfaces of bounded negative curvature always have growth.

\end{abstract}

\section*{Introduction}

Let $G$ be a connected graph and $o\in V(G)$ be a vertex. Denoting by $B_n(o,G)=B_n(o)$ the ball of radius $n$ around $o$ in $G$, the {\it upper growth rate} of $G$ is defined as $\underset{n\to \infty}{\rm limsup}|B_n(o)|^{1/n}$, and this quantity does not depend on the choice of $o$. We define the {\it lower}) {\it growth rate} similarly, with liminf replacing limsup. When they are equal, that is the limit $\underset{n\to \infty}{\rm lim}|B_n(o)|^{1/n}$ exists, we call it the {\it growth rate} of $(G,o)$.

A bounded-degree random rooted graph $(G,o)$ is called a {\it unimodular random graph} if 
it the ``handshake lemma'' of finite graphs is satisfied in a measurable sense. Examples include Benjamini-Schramm limits of finite graphs, and the cluster of the origin by a group-invariant percolation on a Cayley graph.
See \cite{AL} for a thorough definition and equivalents.   
Vadim Kaimanovich asked if it is true that in any unimodular random rooted graph the growth rate exists almost surely. This was refuted by the second author by a counterexample which is a one ended tree with lower growth rate being $1$ {\cite{T}}.
Abért, Fraczyk and Hayes showed that in a unimodular random rooted tree with degree bounded above by $d$, if the upper growth rate is at least $\sqrt{d-1}$, then the growth rate actually exist {\cite{AFH}}. In the closing paragraph they ask if one can drop the condition that the graph is a tree. They expected the answer to be negative. In this paper we confirm this, and show that even the maximal possible upper growth rate can be achieved, while having no growth rate.

\begin{theorem}\label{t.main}
If $d\geq 3$, then there exists a unimodular random graph $({\cal U},o)$ with degree bounded above by $d$ whose upper growth rate is a.s.\ $d-1$ and whose lower growth rate is a.s.\ $1$. 
\end{theorem}
The existence of such an extreme example may be surprising. On one hand, for a Cayley graph $G$ of degree $d$, the growth rate (which always exists) is strictly smaller than $d-1$ whenever $G$ is not a regular tree. One could think that similarly, if a unimodular random graph is not cycle-free, then a positive density of edges are useless from the point of view of growth, because they are ``closing cycles'', and hence the upper growth should be strictly less than $d-1$. Even more so if the unimodular random graph is hyperfinite, as will be the example for Theorem \ref{t.main}, because hyperfiniteness (an equivalent to the suitably defined version of amenability) really enables one to take averages and talk about ``positive density'' in the sense of averaging over space. As the theorem shows, the above heuristics is wrong.

We will proceed by
starting from a weaker construction, then building a stronger one, and finally prove Theorem \ref{t.main} using these previous two.

\medskip

\noindent
{\bf Construction 1}
If $d\geq 3$, then there exists a unimodular random graph $(\Gamma_d,o)$ with degree bounded above by $d+2$ whose upper growth rate is a.s.\ at least $d-1$ and whose lower growth rate is a.s.\ $1$.

\medskip

\noindent
{\bf Construction 2}
If $d\geq 3,\eps>0$, then there exists a unimodular random graph $({\cal U}_{d,\eps},o)$ with degree bounded above by $d$ whose upper  growth rate is a.s.\ at least $(1-\eps)(d-1)$ and whose lower growth rate is a.s.\ $1$.

\medskip

All the above constructions are hyperfinite, but we will show a modification that gives a non-hyperfinite example:

\begin{theorem}\label{t.main_2}
Let $d\geq 6$. There exists a non-hyperfinite unimodular random graph with maximal degree $d$ which has no growth rate and whose upper growth rate is at least $d-4$.
\end{theorem}
This result may be of interest because of its connection to a conjecture of Mikl\'os Ab\'ert (personal communication), that the growth of a unimodular Riemannian surface of bounded negative curvature always exists. Note that one can always assign a unimodular random non-hyperfinite graph to such a surface through an invariant net. Refer to \cite{B} for an example where a graph construction is used to build a Riemannian surface that serves as a counterexample just as the graph does. However, Benjamini's method does not seem to directly apply here.

To ensure that our examples are unimodular, we will rely on the following:
If $(U,o)$ is an unimodular random graph and $(W,o)$ is a random rooted graph on $V(U)$ (where the distribution of $o$ is the same in the two graphs) with the property that $W$'s edge set can be obtained from $U$ by applying some rules of local modifications (such as a factor of iid), then $(W,o)$ is unimodular as well. The reason is that if a mass transport contradicting the mass transport principle existed on $(W,o)$, then that already existed on $(U,o)$.

All of our graphs are going to be simple.
For a graph $G$, let $\gamma(G)$ be its girth (which is the length of the shortest cycle of $G$). If $d$ and $n$ are both odd, then there cannot exists a $d$-regular graph on $n$ vertices; however, with some abuse of terminology, if all but one vertex of a graph have degree $d$ and the exceptional vertex has degree $d-1$, then we still will say that the graph is $d$-regular.
If $G$ is $d$-regular and $r\leq \frac{\gamma(G)}{2}$ is an integer, 
then $|B_r(v,G)|\geq(d-1)^r$, that is $|B_r(v,G)|^{1/r}\leq (d-1)$ for all $v\in V(G)$.
Connected $d$-regular graphs on $n$ vertices with girth at least $(1-o(n)){\rm log}_{d-1}(n)$ are known to exists, see \cite{LS} for a recent simple construction.
In what follows we treat $d$ as fixed and when we mention the ${\rm log}$ function we will mean it in base $(d-1)$. Throughout the paper ${\cal T}$ will denote the {\it canopy tree}, defined as the infinite tree where vertices can be partitioned into generations $0,1,\ldots$ such that generation 0 contains all the leaves, each having a neighbor in generation 1, and every vertex of generation $k\geq 1$ has exactly 2 neighbors in generation $k-1$, and one neighbor in generation $k+1$. Taking the root $o$ to be a vertex of generation $k$ with probability $2^{-k}$ makes $(T,o)$ a unimodular random graph.

\section{Construction 1}

For a vertex $v$ let $ {P_{v}^{\infty}}$ be the unique infinite path of ${\cal T}$ which starts at $v$, let $\T(v)$ be the set of vertices $u$ for which $P_u^{\infty}$ contains $v$, let $v$'s {\it index} ${\tt ind}(v)$ be the shortest distance between $v$ and a leaf of ${\cal T}$ and for 
an edge $\{u_1,u_2\}=e \in E({\cal T})$ let its index ${\tt ind}(e)$ be ${\rm min}\{{\tt ind}(u_1),{\tt ind}(u_2)\}$.

Let $l_0<l_1<\dots $ be a sequence of positive integers such that $l_0\geq d$ and 
$$l_{i+1}-l_i\geq 
2^{2^{l_{i}+2}}{\rm log}_{2}(d-1).$$
Let $I\subset E({\cal T})$ be the set of edges for whose index is contained in $\{l_0,l_1,\ldots\}$. 

If we remove the edges of $I$ from ${\cal T}$, then we get a disconnected graph consisting of infinitely many finite components (we will often call them {\it clusters}). Let ${\tt comp}(I)$ denote the set of these clusters.
For any cluster $K \in {\tt comp}(I)$ there is a unique vertex ${\tt top}(K) \in V(K)$ whose index is maximal in $V(K)$. 
Let ${\tt leaf}(K)$ be the set of leaves of $K$.

If $K\in {\tt comp}(I)$, then ${\tt top}(K)$ is incident to an edge $e\in I$ and any $v\in {\tt leaf}(K)$ is incident to some $f\in I$ (unless $K$ is ``at the bottom", that is ${\tt ind}(e)=l_0$), such that ${\tt ind}(e)=l_{j+1},{\tt ind}(f)=l_j$, for some $j$. Then $|V(K)|=2^{l_{j+1}-l_j}-1$ and $|\T(v)|=2^{l_j+2}-1$.

By the condition on the growth on the sequence $(l_i)$ we have that ${\rm log}_{d-1}|V(K)|\geq 2^{2^{l_j+2}}$, 
thus 
$$\frac{|\T(v)|}{{\rm log}|V(K)|}\leq \frac{2^{l_j+2}-1}{2^{2^{l_j+2}}}.$$
Since $l_j\geq j$, this implies that for any $K\in {\tt comp}(I)$ and leaf $v$ in it,
\begin{equation}\label{eq1}
\frac{|\T(v)|}{{\rm log}|V(K)|}\to 0
\end{equation}
as $j\to \infty$.

For each $K\in {\tt comp}(I)$ we will choose one of two possible {\it types}: ${\tt path}$ or ${\tt exp}$.
Then we will
introduce a new random graph ${\tt new}(K)$ on the vertex set of $K$, depending on the type:

{\it (i),} If $K$ is of type ${\tt path}$, then ${\tt new}(K)$ is a path ${\tt P}_K$ which contains every vertex of $K$ with the following constraints: the leaves of $K$ are visited first (in any order), then $P_K$ visits the remaining vertices of $K$ in any order with the restriction that its last vertex is ${\tt top}(K)$. 

{\it (ii),} If $K$ is of type ${\tt exp}$, then ${\tt new}(K)$ is a $d$-regular graph ${\tt H}_K$ of maximal girth with vertex set $V(K)$.

By the type of a vertex $v\in V({\cal T})$, we will mean the type of the cluster $K$ that contains $v$. 
For any $K\in {\tt comp}(I)$ let the type of ${\tt new}(K)$ be decided by a Bernoulli(1/2), independently from one another.  
Let $({\cal T},o)$ the random rooted graph with the distribution of the root which makes it unimodular.
Define a unimodular random graph ${\cal W}(I)$ on vertex set $V({\cal T})$, with root $o$ and edge set
$$E({\cal W}(I)):=I \cup \underset{K \in {\tt comp}(I)}{\bigcup}E({\tt new}(K)).$$ 
Since all the $K$'s are finite, this can be realized by taking a random uniform choice at each $K$, thus ${\cal W}(I)$ is unimodular.

The following useful observation is straightforward.
\begin{lemma}\label{trivi}
For any two vertices $u,v\in V({\cal T})$, the unique path $P$ between them in ${\cal T}$ and any path $Q$ between them in ${\cal W}(I)$ must go through the same set of edges from $I$ and in the same order. This in particular implies that any two path $Q_1,Q_2$ between them in ${\cal W}(I)$ must go through the same set of edges from $I$.
\end{lemma}

First we show that the lower growth rate of ${\cal W}(I)$ is $1$.

Note that ${\tt leaf}(K)\leq\frac{2}{3}|V(K)|$ for all $K\in {\tt comp}(I)$. Let $c:=\frac{2}{3}$.
For $m\in \N$, let $P$ be a path of minimal length in ${\cal W}(I)$ with the following constraints: the endpoints of $P$ are $o$ and $t={\tt top}(K_t)$ for some $K_t\in {\tt comp}(I)$ such that the type of $K_t$ is ${\tt path}$ and $P$ contains at least $m$ edges $e_1,\dots,e_m$ from $I$ (in this order as moving from $o$ to $v_m$), and let ${\tt ind}(e_1)<{\tt ind}(e_2)<\ldots< {\tt ind}(e_m)$. Let $L$ be the length of $P$.
Observe that $L\geq (1-c) |V(K_t)|$ by the condition that in the path ${\tt P}_{K_t}$ all vertices of ${\tt leaf}(K_t)$ are located in its first segment and $P({o,m})$ must first enter the vertices of $K_t$ in ${\tt leaf}(K_t)$ before it can reach $t$ and must cross all vertices of $V(K_t)\setminus {\tt leaf}(K_t)$ before that, and $|V(K_t)\setminus {\tt leaf}(K_t)|\geq
(1-c)|V(K_t)|$.

We can compare the size of $V(K_t)$ with that of $\T(t)$ by noting that any vertex of $\T(t)$ which is not in $V(K_t)$ is contained in $\T(u)$ for some $u\in {\tt leaf}(K_t)$.
Thus $$|\T(t)|\leq |V(K_t)|+c\cdot |V(K_t)|\cdot |\T(u)|,$$ where $u$ is an arbitrary vertex from ${\tt leaf}(K_t)$.

Because of Lemma \ref{trivi} and the fact that $P$ is minimal, we have $B_L(o,{\cal W}(I)) \subset \T(t)$. By the growth condition on the clusters we have that $|\T(u)|\leq |V(K_t)|$. Combining the above we get that $$|B_L(o,{\cal W}(I))|\leq |\T(t)| \leq |V(K_t)|+c\cdot|V(K_t)|\cdot |V(K_t)|\leq 3|V(K_t)|^2\leq 3\left(\frac{1}{1-c} L\right)^2,$$ and hence $|B_L(o,{\cal W}(I))|^{1/L}\leq \left(3 \left(\frac{ L}{1-c}\right)^2\right)^{1/L}.$
We can let $L\to \infty$ by letting $m\to \infty$ and then this expression goes to $1$.

Next we show that the upper growth rate is high.
For $m\in \N$, let $P$ be a minimal path (in ${\cal W}(I)$) with the following constrains: $P$ is between $o$ and a vertex $\ell$ such that $\ell$ is a leaf in the $K_\ell\in {\tt comp}(I)$ that contains it, 
where the type of $K_\ell$ is ${\tt exp}$ and $P$ contains at least $m$ edges of $I$ of increasing indices. Let the length of $P$ be $L$.

Since all vertices of $P$ are contained in $\T(\ell)$, \eqref{eq1} 
implies that $L\leq \eps_m r$, where $r=\frac{{{\rm log}|V(K_t)|}}{2}$ for some $\eps_m\to 0 $, as $m\to \infty$. 

We also have $$|B_{L+r}(o,{\cal W}(I))|\geq |B_r(u_{o,m},{\tt new}(K))| \geq (d-1)^r.$$
These imply that $$|B_{L+r}(o,{\cal W}(I))|^{\frac{1}{(L+r)}}\geq (d-1)^{\frac{r}{L+r}}\geq (d-1)^{\frac{1}{1+\eps_m}},$$ where $\eps_m\to 0$ as $m\to \infty$. Thus the upper growth rate of ${\cal W}(I)$ is at least $(d-1)$.

Note that the maximal degree of ${\cal W}(I)$ is $d+2$, since in case the type of ${\tt new}(K)$ is ${\tt exp}$, the leaves of $K$ have degree $d$ within ${\tt new}(K)$ (with the potential exception of a single vertex of degree $(d-1)$ if forced by parity constraints) and they will have two extra edges of $I$ incident to them (unless $K$ is at ``the very bottom" of ${\cal T}$).

\section{Construction 2}

Let us call an edge set $J\subset E({\cal T})$ {\it unbounded} if for every $v\in V({\cal T})$, the path $P^{\infty}
_v$ contains at least one edge from $J$ (since this condition must hold for all vertex $v$, this path must contain infinitely many edges from $J$).
$I$ in the previous section was unbounded. 
If $J\subset I$ is unbounded, it will still hold that removing $J$ will result in a graph which consists of infinitely many finite clusters (we will also denote the set of these clusters by ${\tt comp}(J)$).

For $K\in {\tt comp}(J)$ let ${\tt out}_J(K)$ be the set of those vertices $v\in V(K)$ which are different from ${\tt top}(K)$ and which are incident to an edge of $J$.
If we replace $I$ by an unbounded  $J\subset I$, and ${\tt leaf}(K)$ by ${\tt out}_J(K)$ for $K\in {\tt comp}(J)$, then we could repeat the definition of ${\cal W}(I)$ verbatim and would get another unimodular graph ${\cal W}(J)$.

For a vertex $u\in V({\cal T})$, let $K_I(u)$ and $K_J(u)$ denote the clusters containing $u$ in ${\tt comp}(I)$ and ${\tt comp}(J)$ respectively. 
By $J\subset I$ we have $V(K_I(u))\subset V(K_J(u))$.
Moreover, if $v_1\in {\tt leaf}(K_I(u))$ and $v_2\in {\tt out}_J(K_J(u))$, then ${\tt ind}(v_1)\geq {\tt ind}(v_2)$, thus $|\T(v_1)|\geq |\T(v_2)|$.
This implies that $$\frac{\T(v_2)}{{\rm log}|V(K_J(u))|}\leq \frac{\T(v_1)}{{\rm log}|V(K_I(u))|},$$
and hence that \eqref{eq1} can only get stronger, and hence the upper growth rate does not decrease when we move from $I$ to $J$. 
In the previous construction we also used that the fraction of leaves in a cluster $K$ is not greater than $c\cdot |V(K)|$ for some $c<1$, for the proof that the lower growth rate is $1$. In the case of ${\cal W}(J)$ we need a similar property when the leaves of $K$ are replaced by ${\tt out}(K)$. Notice that this also holds when we move from $I$ to $J$ with a potentially even better (that is, smaller) $c<1$.

Given $\eps\in \left(0,\frac{1}{d(d-1)^2}\right)$, we will show that there exists an unbounded $J\subset I$ such that:
\begin{enumerate}[label=(\arabic*)]
    \item For any $v\in V({\cal T})$ there is at most one edge $e_v$ in $J$ which is incident to $v$.
    \item For any cluster $K\in {\tt comp}(J)$ we have
$$|{\tt out}(K)|\leq \frac{\eps}{2} |V(K)|.$$
\end{enumerate}

Suppose first that we found such a $J$. We introduce a new graph ${\tt new}^{*}(K)$ on $V(K)$ for each $K\in {\tt comp}(J)$. The types of $K$ can again be ${\tt path}$ or ${\tt exp}$ (and this decision is made by a Bernoulli(1/2) independently for each clusters).
In case $K$ is of type ${\tt path}$, we use the same definition as before to set ${\tt new}^{*}(K)$ (with ${\tt out}_J(K)$ replacing ${\tt leaf}_J(K)$).
In case $K$ is chosen to be of type ${\tt exp}$, it is defined as follows. Let us select uniformly and random a distinguished subset ${\tt ext}(K)\subset (K\setminus {\tt out}(K))$ of size $\frac{\eps}{2}|V(K)|$. After this, we place a maximal girth $d$-regular graph ${\tt H}_K$ on the vertex set $V(K)$ with the restriction that any two
vertices of ${\tt out}_J(K)\cup {\tt ext}(K)$ are at distance at least $3$ from one another in ${\tt H}_K$ (this is where we need $\eps<\frac{1}{d(d-1)^2}$).
Recall that we may allow a single exceptional vertex of degree $(d-1)$ in case parity constraints force it. When this is the case, then choose the exceptional vertex in such a way that its distance from ${\tt out}_J(K)\cup {\tt ext}(K)$ is at least $2$.
For every vertex $v\in {\tt out}_J(K)\cup {\tt ext}(K)$, let us delete a uniformly chosen edge $e_v$ from ${\tt H}_K$ which is incident to $v$.
Let us call the resulting random graph ${\tt H}_K^{\tt GW}$, and define ${\tt new}^{*}(K)$ as ${\tt H}_K$ if the type of $K$ is ${\tt exp}$.

\begin{remark} The set ${\tt ext}(K)$ will only be important in the next section, for the strongest construction with upper growth rate $(d-1)$, but it is convenient to define it here.
\end{remark}

This concludes the definition of ${\tt new}^{*}(K)$. ${\tt H}_K^{\tt GW}$ is not $d$-regular, but for $r=\frac{\gamma({\tt H}_K)}{2}$, $|B_{r}(v,{\tt H}_K^{\tt GW})|$ stochastically dominates the size of the  first $r$ generation of a supercritical Galton-Watson tree in which every individual has $(d-1)$ children with probability $(1-\eta)$ and $(d-2)$ children with probability $\eta$ (with an appropriate $\eta$ which goes to zero as $\eps\to 0$). Note that this Galton-Watson tree always survives, and hence by the martingale convergence theorem there exists a random variable $W$ such that $\underset{n\to \infty}{\rm lim}\frac{Z_n}{\mu^n}=W$ holds almost surely, with $W>0$. (Here $Z_n$ is the $n$'th generation of the process and $\mu$ is the expected number of children. See e.g. Corollary 5.7 in \cite{LPbook}.) 
This implies that there is an $\eps^{*}$ tending to 0 as $\eps\to 0$, such that for any $v\in {\tt out}_J(K)$, with high probability
$$|B_r(v,{\tt H}_K^{\tt GW})|\geq ((1-\eps^{*})(d-1))^{r}.$$

If we define ${\cal W}^{*}(J)$ as the graph with vertex set $V({\cal T})$ and edge set  
\begin{equation}\label{Wstar}
E({\cal W}^{*}(J)):=J \cup \underset{K \in {\tt comp}(J)}{\bigcup}E({\tt new}^{*}(K)),
\end{equation}
then essentially the same argument as in the previous section shows that it satisfies the properties claimed for Construction 2 in the Introduction.
A minor technical difference is that $|B_r(v,{\tt H}_K^{\tt GW})|\geq ((1-\eps^{*})(d-1))^{r}$ only holds with high probability, but a Borel-Cantelli argument works to show that as we take the analogues of the paths $P$ we used before, we will encounter infinitely many vertices $v$ for which this holds, which is enough to reach the conclusion that the upper growth rate is $(1-\eps^{*})(d-1)$.

To find a $J$ as above, let us do the following.
Say that two vertices $u,v$ of the same index are $L$-equivalent if $P_u^{\infty}$ and $P_v^{\infty}$ coincide no later than after $L$ steps. 
In this way every equivalence class is finite.
Select $L$ to be greater than $1/{\eps}$ and let $k$ be large enough such that $l_{k+1}-l_k>L$ ($l_i$ is the sequence we used for defining $I$). For all $l_m$ for which $m>k$, select uniformly at random a representative $u$ from each equivalence class of vertices whose index is $l_m$, and put the first edge of $P^{\infty}_u$ into $J$. 
It is easy to see that $J$ satisfies the required conditions (we may use a Borel-Cantelli argument to show that $J$ is unbounded a.s.). We highlight the fact that in ${\cal T}\setminus J$ (and hence in ${\cal W}^{*}(J)\setminus J$) all clusters reach the leaf-level of the canopy tree, and hence it is not surprising that only an ``$\eps$ fraction of the vertices" are connected to other clusters.

\section{Construction with the maximal possible upper growth rate}

In this section we prove Theorem~\ref{t.main}. 
Recall the notion of Cartesian product $G\square H$ of two graphs $G$ and $K$.
The vertex set of it is $V(G\square H)=V(G)  \times V(H)$ and a vertex $(g_1,h_1)$ is adjacent to $(g_2,h_2)$ if
either $g_1=g_2$ and $h_1$ is connected to $h_2$, in which case we will say that the edge $\{(g_1,h_1),(g_2,h_2)\}$ is {\it vertical},
or $g_1$ is adjacent to $g_2$ and $h_1=h_2$, in which case we say that the edge $\{(g_1,h_1),(g_2,h_2)\}$ is {\it horizontal}.
We first consider ${\cal T} \square {\cal T}$, where ${\cal T}$ is the canopy tree.

For $v\in V({\cal T})$, let $v^+$ be the next vertex after $v$ on $P^{\infty}_v$.
For $(v,u)={\mathbf v}\in V({\cal T}\square {\cal T})$, let ${\mathbf v}_{+}=(v^+,u)$, let $f^\uparrow({\mathbf v})$ be the vertical edge $\{(v,u),(v,u^+)\}\in E({\cal T}\square {\cal T})$ and $f^\rightarrow({\mathbf v})$ be the horizontal edge $\{(v,u),(v^+,u)\}=\{{\mathbf v},{\mathbf v}_+\} \in E({\cal T}\square {\cal T})$.
Let us say that a path $P$ starting at ${\mathbf v}\in V({\cal T}\square {\cal T})$ and ending at ${\mathbf u}$ is {\it upward} if every step is made toward infinity in either the horizontal or in the vertical direction, that is, for any ${\mathbf x}\in V(P)$ the edge $e\in E(P)$ which separates ${\mathbf x}$ from ${\mathbf u}$ is either $f^\uparrow({\mathbf x})$ or $f^\rightarrow({\mathbf x})$.
For ${\mathbf u}\in V({\cal T}\square {\cal T})$, let us define the following vertex set:
$${\tt t}_{\square}({\mathbf u}):=\{{\mathbf v}\in V({\cal T}\square {\cal T})|\text{there exists an upward path from ${\mathbf v}$ to ${\mathbf u}$}\}.$$

For $v\in V({\cal T})$ let ${\cal T}_v$ be the subgraph (or {\it fiber}) of ${\cal T}\square {\cal T}$ induced by $\{(v,u): u \in V({\cal T})\}$. 
Suppose that for every such $v$ an unbounded edge set $J_{v}\subset E({\cal T}_{v})$ is given. For an arbitrary $(v,u)={\mathbf v}\in V({\cal T}\square {\cal T})$, let ${\tt comp}_{J_{v}}({\mathbf v})$ denote the vertex set of the component $K\in {\tt comp}({J_{v}})$, for which ${\mathbf v}\in V(K)$.
For $f\in E({\cal T}_{v})$, where $f=f^\uparrow({\mathbf u})$,
let us say that {\it $f$ is turnable with respect to $J_{v^+}$}, if:
\begin{equation}\label{eq1new} 
|{\tt t}_{\square}({\mathbf u}_+)|\leq \eps(v^+){\rm log}(|{\tt comp}_{J_{v^+}}({\mathbf u}_+)|).
\end{equation}
The role of this equation will be similar to that of \eqref{eq1} in the previous constructions.

Set $\eps_k:=\frac{1}{4^kd(d-1)^2}$ and for a vertex $v\in V({\cal T})$ let $\eps(v)=\eps_{{\tt ind}(v)}$.
For every fiber ${\cal T}_v$ consider the $J_v$ used in Construction 2. So $J_v$ is an unbounded edge set within ${\cal T}_v$ and we constructed ${\cal U}_{\eps(v)}$ as ${\cal W}^{*}(J_v)$ in \eqref{Wstar}. Note that if we replace $J_v$ with an unbounded subset of it, then it will still satisfy (1) and (2) in Construction 2, and hence the new ${\cal W}^{*}(J_v)$ is an example with the properties of Construction 2. Our first aim is to find such a replacement in such a way that we also have many turnable edges.

Observe that for any $f\in E({\cal T}_v)$, there exists an $N(f)\in \N$ (depending only on ${\tt ind}(f)$ and ${\tt ind}(v)$), such that if $J_{v^+}$ only contains edges whose index is either strictly less than ${\tt ind}(f)$ or strictly greater than $N(f)$, then $f$ is turnable with respect to $J_{v^+}$.
Using this observation we can find (using some local algorithm) a family of unbounded edge sets within $J'_v\subset J_v$, such that for each $v\in V({\cal T})$ the set $J_v^{\tt turn}\subset J'_v$ of edges that are turnable with respect to $J'_{v^+}$ is unbounded. We can define the $J'_v$ by induction on ${\tt ind}(v)$ and choosing $J'_u$ to be $J_u$ for $u$ with ${\tt ind}(u)=0$.
Having defined the unbounded edge sets $J'_v\subset E({\cal T}_v)$ for all $v\in V({\cal T})$, we have the graph ${\cal W}^{*}(J'_v)$ on the vertex set $V({\cal T}_v)$ with its distinguished subset 
$\bigcup_{K\in {\tt comp}(J'_v)}({\tt ext}(K))$ of vertices of degree at most $(d-1)$, for all $v\in V({\cal T})$.

Let us call a vertex ${\mathbf v}\in {\cal T}_u$ {\it lucky}, if $f^\uparrow(\mathbf v)\in J_u^{\tt turn}$, ${\mathbf v}_+\in \bigcup_{K\in {\tt comp}(J_{{u^+}}')}({\tt ext}(K))$, and for the other vertex ${\mathbf v}'\in {\cal T}_{u'}$ for which $({\mathbf v}')_+={\mathbf v}_+$, we have that $f^\uparrow({\mathbf v}')\not \in J'_{u}$.

Now we are ready to define ${\cal U}$ and prove Theorem~\ref{t.main}. The vertex set is $V({\cal U})=V({\cal T}\square {\cal T})$ with the distribution of the root ${\mathbf o}$ which makes $({\cal T}\square {\cal T},{\mathbf o})$ unimodular.
To define the edge set of ${\cal U}$, let ${E^\uparrow}=\cup_{v\in V({\cal T})} E({\cal W}^{*} (J'_v))$.
If ${\mathbf v}$ is lucky, we replace $f^\uparrow({\mathbf v})\in E^\uparrow$ by $f^\rightarrow({\mathbf v})$. Formally, if 
$L\subset V({\cal T}\square {\cal T})$ is the set of lucky vertices, $F^\rightarrow=\bigcup_{{\mathbf v}\in L}f^\rightarrow({\mathbf v})$ and $F^\uparrow=\bigcup_{{\mathbf v}\in L}f^\uparrow(\mathbf v)$, then define $E({\cal U})$ to be $F^\rightarrow\cup (E^\uparrow\setminus F^\uparrow)$.
 
Observe that if $P$ and $Q$ are paths in ${\cal U}$ from ${\mathbf v}$ to ${\mathbf u}$, then they must go through the same set of edges from $F^\rightarrow$ and also the same set of edges from $J'_v$ for all $v$.

To show that ${\cal U}$ has the desired properties from Theorem~\ref{t.main}, note that for ${\mathbf v}\in V({\cal T}_u)$, conditioned on $\lambda({\mathbf v})\in J_v^{\tt turn}$, the probability that ${\mathbf v}$ is lucky is $\frac{\eps (u)}{2}$, a constant depending only on ${\tt ind}(u)$, independently for vertices within ${\cal T}_u$. This implies that any infinite path within ${\cal T}_u$ must encounter a lucky vertex which is connected to ${\cal T}_{u^+}$ in ${\cal U}$.

After this, essentially the same proof as in the case of Construction 1 works (with \eqref{eq1} replaced by \eqref{eq1new}), thus we only sketch it.
To show that the upper growth rate is high, we consider a path $P_{{\mathbf o},m}$ of minimal length $L=L(m)$ with the following properties:
$P_{{\mathbf o},m}$ starts at ${\mathbf o}$, crosses at least $m$ edges from $F^\rightarrow$, its last edge is in $F^\rightarrow$, and ends at a vertex ${\mathbf x}$ in a cluster $K\in {\tt comp}(J'_w)$ of some $J'_w$ of type ${\tt exp}$. With high probability the ball $B_r({\mathbf x},{\cal U})$ of radius $r=\frac{{\rm log}(|V(K)|)}{2}$ has at least $((1-\eps(w))(d-1))^r$ many vertices.
By our construction the last edge on $P_{{\mathbf o},m}$ is in ${\cal U}$ 
so the vertex ${\mathbf y}$ right before ${\mathbf x}$ had $f^\uparrow ({\mathbf y})$ turnable (and so $f^\uparrow ({\mathbf y})$ is not in $P_{{\mathbf o},m}$ but $f^\rightarrow ({\mathbf y})$ is). Thus all of $P_{{\mathbf o},m}$ is contained ${\tt t}_{\square}({\mathbf x})$ and since $f^\uparrow ({\mathbf y})$ was turnable, we have that $\frac{L}{r}\to 0$ as $m\to \infty$, thus $|B_{r+L}({\mathbf x},{\cal U})|^{\frac{1}{r+L}}\geq ((1-\eps(w))(d-1))^{\frac{r}{r+L}}$. So ${\rm limsup}_{m\to \infty}|B_{r+L}({\mathbf x},{\cal U})|^{\frac{1}{r+L}}\geq (1-\eps(w))(d-1)$, but this is true for any $w\in V({\cal T})$ so we got that the upper growth rate is $(d-1)$.

To show that the lower growth rate is $1$, we consider a path $Q_{{\mathbf o},m}$ of minimal length $l=l(m)$ with respect to the following properties:
$Q_{{\mathbf o},m}$ starts at ${\mathbf o}$, crosses at least $m$ edges from $F^\rightarrow$ and ends at a vertex ${\mathbf x}={\tt top}(K)$ for a cluster $K\in {\tt comp}(J'_w)$ of type ${\tt path}$. We have that $l\geq (1-\eps(w))|V(K)|$. By noting that $B_l({\mathbf o},{\cal U})\subset V(K)\cup\underset{{{\mathbf u}:{\mathbf u}_+\in V(K)}}\bigcup{\tt t}_{\square}({\mathbf u})$, we can also estimate $|B_l({\mathbf o},{\cal U})|\leq |V(K)|+|V(K)|\cdot {\rm log}(|V(K)|)\leq 2\cdot \left(\frac{l}{1-\eps(w)}\right)\cdot {\rm log}\left(\frac{l}{1-\eps(w)}\right)$. This implies that $|B_l({\mathbf o},{\cal U})|^{\frac{1}{l}}\to 1$ as $m\to \infty$.

\section{A non-hyperfinite example}

We prove Theorem \ref{t.main_2} by showing that the
Cartesian product ${\cal U}\square \T_3$ (where $\T_3$ is the $3$-regular tree) also has no growth rate and its upper growth rate is the same as that of ${\cal U}$ for $d\geq 4$, while its maximal degree is the maximal degree of ${\cal U}$ plus $2$ (and it is not-hyperfinite since it contains non-hyperfinite subgraphs). Let $o=(o_1,o_2)$ be the root in the Cartesian product.

The upper growth rate of ${\cal U}\square \T_3$ is clearly at least that of ${\cal U}$.
To get an estimate on the lower growth rate let us use the fact that any vertex of distance $r$ can be reached by first moving $i$ steps horizontally and then at most $r-i$ steps vertically.
This gives us the following estimate:

$$|B_r(o,{\cal U}\square \T_3)|\leq |B_r(o_1,{\cal U})|+3\sum_{i=1}^{r}2^{i-1}|B_{r-i}(o_1,{\cal U})|$$

Using the fact that the lower growth rate of ${\cal U}$ is $1$, for any $c>0$ there is a sequence $r_1,r_2,\dots$, where $r_i\to \infty$ such that $|B_{r_i}(o_1)|<(1+c)^{r_i}$. Fix such a $c$ such that $2+2c$ is strictly less than the upper growth rate of ${\cal U}$.

A very rough upper estimate on $|B_r(o,{\cal U}\square \T_3)|$ will be sufficient (we simply replace $2i$ by $2^r$ and $B_{r-i}(\pi_1(o),{\cal U})$ by $B_{r}(\pi_1(o),{\cal U})$):
$$|B_r(\pi_1(o),{\cal U})|+3\sum_{i=1}^{r}2^{i-1}|B_{r-i}(\pi_1(o),{\cal U})|\leq $$
$$\leq  |B_r(\pi_1(o),{\cal U})|+3\sum_{i=1}^{r}2^{r}|B_{r}(\pi_1(o),{\cal U})|=$$
$$=(1+3r2^r)|B_r(\pi_1(o),{\cal U})|$$
$$\leq 6r2^r|B_r(\pi_1(o),{\cal U})|.$$

If we use this estimate for some $r=r_i$ from our sequence, we know that $|B_r(o_1,{\cal U})|<(1+c)^r$, so we have:
$$|B_r(o,{\cal U}\square \T_3)|\leq 6r2^r|B_r(\pi_1(o),{\cal U})|\leq 6r(2(1+c))^r.$$
Since the upper growth rate of ${\cal U}$ is strictly greater than $2+2c$, we have that the upper growth rate of ${\cal U}\square \T_3$ is strictly greater than its lower growth rate.
Using our stronger construction we obtain a non-hyperfinite graph of no growth with maximal degree $d\geq 6$ whose upper growth rate is arbitrarily close to $d-4$.

\bigskip

\noindent
{\bf Acknowledgements.} This research was supported by the Icelandic Research Fund grant No. 239736-051 and ERC
grant No. 810115-DYNASNET.

\ \\
{\bf P\'eter Mester}\\
Division of Mathematics, The Science Institute, University of Iceland\\
Dunhaga 3 IS-107 Reykjavik, Iceland\\
\texttt{pmester[at]alumni.iu.edu}\\

\noindent
{\bf \'Ad\'am Tim\'ar}\\
Division of Mathematics, The Science Institute, University of Iceland\\
Dunhaga 3 IS-107 Reykjavik, Iceland\\
and\\
HUN-REN Alfr\'ed R\'enyi Institute of Mathematics\\
Re\'altanoda u. 13-15, Budapest 1053 Hungary\\
\texttt{madaramit[at]gmail.com}\\

\end{document}